%% file: Invariants_and_Gorenstein_projective_modules.tex
\UseRawInputEncoding
\documentclass[oneside, a4paper,reqno]{amsart}
\usepackage{pdfsync}
\usepackage{stmaryrd}
\usepackage{mathrsfs}
\usepackage{multicol}
\usepackage{amsmath, amsthm, amscd, amssymb, latexsym, eucal}
\usepackage[all]{xy}
\def\serieslogo@{} \def\@setcopyright{} \makeatother
\usepackage{multienum}

\usepackage{hyperref}
\usepackage{color}
\usepackage{cite}

\hypersetup{colorlinks=true,
     breaklinks=true,
     linkcolor=blue,    % how to specify linkcolor according to chapter / contents ?
     bookmarks=true,
      pageanchor=true}

\makeatletter
\renewcommand*\env@matrix[1][c]{\hskip -\arraycolsep
  \let\@ifnextchar\new@ifnextchar
  \array{*\c@MaxMatrixCols #1}}
\makeatother

\usepackage{color}

\input{rgb}

\usepackage[colorinlistoftodos]{todonotes}

\numberwithin{equation}{section}
\newtheorem{thm}{Theorem}[section]

\newtheorem{lem}[thm]{Lemma}
\newtheorem{prop}[thm]{Proposition}

\theoremstyle{definition}
\newtheorem{defn}[thm]{Definition}

\newtheorem{exam}[thm]{Example}

\newcommand{\lxr}{\longrightarrow}
\newcommand{\iso}{\cong}
\newcommand{\aGP}{\mathscr GP}
\DeclareMathOperator*{\Ker}{\mathsf{Ker}}

\newsavebox{\proofbox}
\savebox{\proofbox}{\begin{picture}(7,7)%
  \put(0,0){\framebox(7,7){}}\end{picture}}

\begin{document}
%\linenumbers

\title[]{Invariants and Gorenstein projective modules}
\author[]{Nan Gao$^{*}$, Chi-Heng Zhang
}
\address{Department of Mathematics, Shanghai University, Shanghai 200444, PR China}
\thanks{* is the corresponding author.}
\email{nangao@shu.edu.cn, zchmath@shu.edu.cn }
\thanks{Supported by the National Natural Science Foundation of China (Grant No. 11771272 and 11871326).}

\date{\today}



\begin{abstract}\ Invariants with respect to recollements of the stable category of Gorenstein projective $A\mbox{-}$modules over an algebra $A$ and stable equivalences are investigated. Specifically, the Gorenstein rigidity dimension is introduced. It is shown that the Gorenstein rigidity dimension is invariant with respect to both Morita equivalences and the stable equivalences of Gorenstein projective modules. As a consequence,
the Gorenstein rigidity dimension is shown the invariant of derived equivalences. The Gorenstein rigidity dimension is compared along the recollements of the stable category of Gorenstein projective modules. Moreover, the bounds of Gorenstein rigidity dimension is given
for several classes of algebras, respectively.
\end{abstract}

\maketitle

%\setprotcode\font
%    {\it \setprotcode\font}
%    {\bf \setprotcode\font}
%    {\bf \it \setprotcode\font}
%    \pdfprotrudechars=1

\vskip 20pt

\section{\bf Introduction}

\vskip 10pt

For Gorenstein homological algebra we refer to \cite{EJ1, EJ2, GZ, H}, whose main idea is to replace projective modules by Gorenstein
projective modules. These modules were introduced by Enochs and Jenda \cite{EJ1} as a generalization of finitely generated module of $G$-dimension zero over a two-sided noetherian ring, in the sense of Auslander and Bridger \cite{AB}. It is known that an algebra has many invariants, for example, the finiteness of global dimension, the finitistic dimension,
and CM-finiteness (see \cite{BR, G, Be}).

\vskip 10pt

It is a well-known result of Rickard \cite{R} which says that a derived equivalence between two selfinjective algebras always induces a stable equivalence of Morita type. So, many basic questions on this aspect arise.  For example, Hu and Pan \cite{HP} introduced  nonnegative and uniformly bounded functors between derived categories of abelian categories, respectively. They proved that a derived equivalence induces a stable equivalence of Gorenstein projective objects. The invariant of stable equivalences of Gorenstein projective modules is of interest in the representation theory of algebras, in Gorenstein homological algebra, and in the theory of singularity categories (see, e.g., \cite{BDZ, BOJ, BR}). How to characterize this is a basic problem.

\vskip 10pt

One can approach these questions by Gorenstein rigidity dimension. One way is by  stable equivalences.
 Another way is by recollements, which have been introduced by Beilinson et al.\cite{BBD}.
A recollement of a triangulated category by another two triangulated categories is a diagram of six functors between
 these categories, generalising Grothendieck's six functors.

\vskip 10pt

In the paper, the notion of Gorenstein rigidity dimension are introduced. Then some properties are collected, which shows that this kind of dimension has some advantages in Gorenstein homological algebra. Recollements and stable equivalences are used as a crucial tool. Suppose that $A$, $B$ and $C$ are Artin algebras. If the stable category $A\mbox{-}\underline{{\rm Gproj}}$ of Gorenstein projective $A$-modules admits a recollement with respect to $B\mbox{-}\underline{{\rm Gproj}}$ and $C\mbox{-}\underline{{\rm Gproj}}$ of $B$ and $C$, then we show that the Gorenstein rigidity dimension of $A$ is bigger than them of $B$ and $C$ (see Theorem~\ref{recollementgridity}). In particular, the Gorenstein rigdity dimension is shown the invariant under both Morita equivalence and the equivalence between $A\mbox{-}\underline{{\rm Gproj}}$ and $B\mbox{-}\underline{{\rm Gproj}}$ (see Proposition~\ref{mainprop} and Theorem~\ref{equivalences}). Based on these, the Gorenstein rigidity dimension of several classes of algebras are described.

\vskip 20pt

\section{\bf Preliminaries}

\vskip 10pt

In this section, we fix notations and recall some basic concepts.

\vskip 10pt

Let $A$ be an Artin algebra. Denote by $A\mbox{-}{\rm Mod}$ (resp. $A\mbox{-}{\rm mod}$) the category of (finitely generated) left $A\mbox{-}$modules, and $A\mbox{-}{\rm proj}$ the full category of finitely generated projective $A\mbox{-}$modules. A module $G$ of $A\mbox{-}{\rm mod}$ is {\bf Gorenstein projective} if there is an exact sequence
$$\cdots \lxr P^{-1}\lxr P^{0}\stackrel{d^{0}}{\lxr} P^{1}\lxr \cdots$$
in $A\mbox{-}{\rm proj}$, which stays exact after applying ${\rm Hom}_{A}(-, P)$ for each module $P$ in $A\mbox{-}{\rm proj}$, such that $G\iso \Ker d^{0}$ (see \cite{EJ1, EJ2}). Denote by $A\mbox{-}{\rm Gproj}$ the full subcategories of Gorenstein projective modules in $A\mbox{-}{\rm mod}$, and $A\mbox{-}\underline{{\rm Gproj}}$ the stable category of $A\mbox{-}{\rm Gproj}$ that modulo $A\mbox{-}{\rm proj}$. Similarly, denote by $A\mbox{-}{\rm Ginj}$ the full subcategory of Gorenstein-injective modules in $A\mbox{-}{\rm mod}$.

\vskip 5pt

Recall from \cite{Be} that an Artin algebra $A$ is of {\bf finite Cohen-Macaulay type} (resp. finite CM-type for simply), if there are only finitely many isomorphism classes of indecomposable finitely generated Gorenstein projective modules. Set $G=\oplus_{i=1}^{n}G_{i}$, where $G_{i}$ are all non-isomorphic indecomposable finitely generated Gorenstein projective modules, and $\aGP(A):=(\rm{End}_{A}(G))^{op}$, which we call $\aGP(A)$ the Cohen-Macaulay Auslander algebra. Recall from \cite{Ch12} that an Artin algebra $A$ is called {\bf CM-free} provided that $A\mbox{-}{\rm Gproj}= A\mbox{-}{\rm proj}$, that is, all its finitely generated Gorenstein projective modules are projective.
Recall from \cite{H} that an Artin algebra $A$ is {\bf Gorenstein} if ${\rm inj.dim} {_{A}A}<\infty$ and ${\rm inj.dim} A_{A}<\infty$.

\vskip 5pt

Now we write $C^{b}(A\mbox{-}{\rm Mod}),\ K^{b}(A\mbox{-}{\rm Mod})$ and $D^{b}(A\mbox{-}{\rm Mod})$ (resp. $C^{b}(A),\ K^{b}(A)$ and $D^{b}(A)$) for the bounded complex category, bounded homotopy category and bounded derived category of $A\mbox{-}{\rm Mod}$ (resp. $A\mbox{-}{\rm mod}$), respectively.

\vskip 5pt

Given triangulated categories $\mathcal{D}, \mathcal{X}$ and $\mathcal{Y}$. Recall from \cite{BBD} that $\mathcal{D}$ is said to be a {\bf recollement} of $\mathcal{X}$ and $\mathcal{Y}$, if there are six triangulated functors as in the following diagram
\[
\xymatrix@C=0.5cm{
\mathcal{Y} \ar[rrr]^{{\rm i}} &&&
\mathcal{D} \ar[rrr]^{{\rm e}}  \ar @/_1.5pc/[lll]_{{\rm q}}  \ar @/^1.5pc/[lll]^{{\rm p}} &&&
\mathcal{X} \ar @/_1.5pc/[lll]_{{\rm l}}  \ar @/^1.5pc/[lll]^{{\rm r}}
 }
\]
such that
\begin{enumerate}
\item \ $({\rm q}, {\rm i}),\ ({\rm i}, {\rm p}),\ ({\rm l}, {\rm e})$ and $({\rm e}, {\rm r})$ are adjoint pairs;

\item \ ${\rm i},\ {\rm l}$ and ${\rm r}$ are fully faithful;

\item \ ${\rm e}{\rm l}=0$;

\item \ For each $Z\in \mathcal{D}$, the counits and units give rise to distinguished triangles:
$${\rm le}(Z)\lxr Z\lxr {\rm iq}(Z)\lxr \ \ \ and\ \ \ {\rm ip}(Z)\lxr Z\lxr {\rm re}(Z)\lxr.$$
\end{enumerate}
A {\bf ladder} of $(\mathcal{Y}, \mathcal{D}, \mathcal{X})$ is a finite or infinite diagram
of triangulated categories and triangle functors:
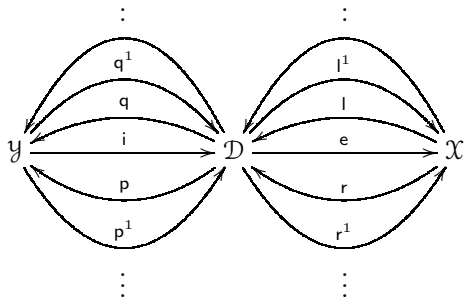
\[
\xymatrix@C=3em@R=3.5em{
& \vdots &  & \vdots & \\
\mathcal{Y} \ar[rr]^{\mathsf{i}} \ar @/^2.3pc/[rr]^{\mathsf{q}^1}   \ar @/_3.0pc/[rr]^{\mathsf{p}^1} &  &
\mathcal{D} \ar[rr]^{\mathsf{e}} \ar @/_1.1pc/[ll]_{\mathsf{q}}     \ar @/^1.5pc/[ll]_{\mathsf{p}}
\ar @/^2.3pc/[rr]^{\mathsf{l}^1}    \ar @/_3.0pc/[rr]^{\mathsf{r}^1}   \ar @/_3.6pc/[ll]_{} &  &
\mathcal{X}  \ar @/_1.1pc/[ll]_{\mathsf{l}}   \ar @/^1.5pc/[ll]_{\mathsf{r}}
\ar @/_3.6pc/[ll]_{}\\
& \vdots &  & \vdots &\\
 }
\]
such that any three consecutive rows form a recollement of triangulated categories. Multiple occurrence of the same recollement is allowed.
The {\bf height} of the ladder $(\mathcal{Y}, \mathcal{D}, \mathcal{X})$ is the number of recollements contained in it.

\vskip 10pt

\section{\bf Gorenstein rigidity dimension}

In the section, we introduce the notion of Gorenstein rigidity dimension. Then we show that it is the invariant
under both Morita equivalences and stable equivalences of Gorenstein projective modules, and it is characterised along the recollements of
stable categories of Gorenstein projective modules. Based on these, the Gorenstein rigidity dimension of several classes of algebras are described.

\vskip 10pt

Before we introduce the notion of Gorenstein rigidity dimension, we need some preparation.

\vskip 5pt

\begin{defn}\ Let $A$ be an Artin algebra and $T$ an $A$-module. The {\bf $T$-dominant dimension} of $M$ in $A\mbox{-}{\rm mod}$, denoted by
$T\mbox{-}{\rm domdim}M$, given by
$$T\mbox{-}{\rm domdim}M={\rm sup} \{n|0 \lxr M \lxr T^{0} \lxr T^{1} \lxr \cdots \lxr T^{n-1} \ is \ an \ exact $$
$$\ sequence\ such \ that \ \ T^{i} \in {\rm add}T \ for \ all \ 0 \leq i \leq n-1\}.$$
\end{defn}

\vskip 10pt

Recall that an $A\mbox{-}$module $M$ is called a generator if $A\in {\rm add}M$, where ${\rm add}M$ is the smallest full subcategory of $A\mbox{-}{\rm mod}$ which contains $M$ and is closed under finite direct sums and direct summands. We say $M$ is a Gorenstein projective generator which means that $M$ is both a generator and a Gorenstein projective $A$-module.

\vskip 10pt

\begin{defn}\ The {\bf Gorenstein rigidity dimension} of an Artin algebra $A$ is defined as follows:
$${\rm Grigdim}A={\rm sup}\{DM\mbox{-}{\rm domdim}\Gamma|M \ is \ a \  Gorenstein \ projective \ generator \ in $$
$$\ A\mbox{-}{\rm mod} \ and \ \Gamma={\rm End}_{A}(M) \}.$$
\end{defn}

\vskip 10pt

Recall from \cite[Definition 2.5]{CFKKY} that the rigidity degree of an $A$-module $M$, denoted by ${\rm rd}_{A}M$, is defined as:
${\rm rd}_{A}M={\rm sup}\{n| {\rm Ext}^{i}_{A}(M, M)=0, \ for \ all \ 1\leq i\leq n\}$. Then we have

\vskip 5pt

\begin{lem}\
\label{domdim}
Let $A$ be an Artin algebra and $M$ a Gorenstein projective generator in $A\mbox{-}{\rm mod}$. Let $\Gamma:={\rm End}_{A}(M)$. Then $DM\mbox{-}{\rm domdim}\Gamma ={\rm rd}_{A}M+2$.
\begin{proof}\
Let
$$0\lxr M\lxr I_{0}\lxr I_{1}\lxr \cdots\lxr I_{n}\lxr \cdots$$
be the minimal injective resolution of $M$. Applying ${\rm Hom}_{A}(M, -)$, we obtain the following complex
$$0\lxr \Gamma\lxr {\rm Hom}_{A}(M, I_{0})\lxr \cdots\lxr {\rm Hom}_{A}(M, I_{n})\lxr \cdots.$$
Since ${\rm Hom}_{A}(M, -)$ is a fully faithful functor and ${\rm Hom}_{A}(M, I_{i})\in {\rm add}DM$, it follows  that $DM\mbox{-}{\rm domdim}\Gamma ={\rm rd}_{A}M+2$.
\end{proof}
\end{lem}

\vskip 10pt

There is a reformulation of the definition, using Lemma~\ref{domdim}:
$${\rm Grigdim}A={\rm sup}\{{\rm rd}_{A}M|M \ is \ a \  Gorenstein \ projective \ generator \ in  \ A\mbox{-}{\rm mod} \}+2.$$
Consequently, ${\rm Grigdim}A\geq 2$ for any algebra $A$.

\vskip 10pt

Recall from \cite{GK} that an Artin algebra $A$ is gendo-Gorenstein if $A$ is isomorphic to the endomorphism algebra of a Gorenstein projective generator over a Gorenstein Artin algebra. Equivalently, there exists an idempotent $e$ of $A$ with
$eAe$ being a Gorenstein algebra and $Ae$ being a Gorenstein projective generator such that $A\cong {\rm End}_{eAe}(Ae)$.

\vskip 10pt

\begin{exam}\
\label{gendo}
Let $(A,e)$ be a gendo-$d$-Gorenstein algebra for some non-negative integer $d$. If ${\rm rd}_{eAe}eA\geq n-2$ for some integer $n\geq 2$, then
$${\rm Grigdim}eAe\geq n.$$
\begin{proof}\ Since $(A,e)$ is a gendo-Gorenstein algebra, we have from Lemma~\ref{domdim} that
$$D(eA)\mbox{-}{\rm domdim}A={\rm rd}_{eAe}eA+2.$$
Since ${\rm rd}_{eAe}eA\geq n-2$, it follows that $D(eA)\mbox{-}{\rm domdim}A\geq n$. By the definition of Gorenstein rigidity dimension,
we have that ${\rm Grigdim}eAe\geq n$.
\end{proof}
\end{exam}

\vskip 10pt

Next we collect some properties of Gorenstein rigidity dimension.

\vskip 5pt

\begin{prop}\
\label{mainprop}
Let $A$ and $B$ be Artin algebras. Then the following hold.

\vskip 5pt
\begin{enumerate}
\item \ ${\rm Grigdim}(A\times B)= {\rm min}\{{\rm Grigdim}A,\ {\rm Grigdim}B\}$.

\vskip 5pt

\item \ If $k$ is perfect, then ${\rm Grigdim}(A\otimes_{k}B)\geq {\rm min}\{{\rm Grigdim}A,\ {\rm Grigdim}B\}$.

\vskip 5pt

\item \ If $A$ and $B$ are Morita equivalent, then ${\rm Grigdim}A={\rm Grigdim}B$.

\end{enumerate}
\begin{proof}\ Let $M$ and $N$ be the Gorenstein projective generators in $A\mbox{-}{\rm mod}$ and $B\mbox{-}{\rm mod}$ respectively. Taking the minimal injective resolution
$$0\lxr M\lxr I^{0}\lxr I^{1}\lxr \cdots\lxr I^{n}\lxr \cdots,$$
of $M$ and applying the functor ${\rm Hom}_{A}(M, -)$, we obtain the following complex
$$0\lxr {\rm End}_{A}(M)\lxr {\rm Hom}_{A}(M, I^{0})\lxr \cdots\lxr {\rm Hom}_{A}(M, I^{n})\lxr \cdots \ \ (*)$$
with ${\rm Hom}_{A}(M, I^{i})\in {\rm add}DM$. Similarly, taking the minimal injective resolution
$$0\lxr N\lxr J^{0}\lxr J^{1}\lxr \cdots\lxr J^{n}\lxr \cdots,$$
of $N$ and applying the functor ${\rm Hom}_{B}(N, -)$, we obtain the following complex
$$0\lxr {\rm End}_{B}(N)\lxr {\rm Hom}_{B}(N, J^{0})\lxr \cdots\lxr {\rm Hom}_{B}(N, J^{n})\lxr \cdots \ \ (**)$$
with ${\rm Hom}_{B}(N, J^{i})\in {\rm add}DN$.

\vskip 5pt

(1) \ Since ${\rm End}_{A\times B}(M\times N)\cong {\rm End}_{A}(M)\times {\rm End}_{B}(N)$ and $D(M\times N)\cong DM\times DN$, we have that $D(M\times N)\mbox{-}{\rm domdim}{\rm End}_{A\times B}(M\times N)= {\rm min}\{DM\mbox{-}{\rm domdim}{\rm End}_{A}(M),$ \ $DN\mbox{-}{\rm domdim}{\rm End}_{B}(N)\}$. This implies that
$${\rm Grigdim}(A\times B)= {\rm min}\{{\rm Grigdim}A,\ {\rm Grigdim}B\}.$$

\vskip 5pt

(2) \ From \cite[Proposition 2.6]{HLXZ}, we know that $M\otimes_{k}N$ is a Gorenstein projective module in $A\otimes_{k}B\mbox{-}{\rm mod}$. Since $A\in {\rm add}M$ and $B\in {\rm add}N$, we have that $A\otimes_{k}B\in {\rm add}(M\otimes_{k}N)$. This implies that $M\otimes_{k}N$ is the Gorenstein projective generator in $A\otimes_{k}B\mbox{-}{\rm mod}$.

\vskip 5pt

Since $k$ is perfect, we have that
$${\rm End}_{A\otimes_{k}B}(M\otimes_{k}N)\cong {\rm End}_{A}(M)\otimes_{k} {\rm End}_{B}(N)$$
as $k\mbox{-}$algebras. Consider the tensor product of complexes $(*)$ and $(**)$. Then by $${\rm Hom}_{A}(M, I^{i})\in {\rm add}DM, \ {\rm Hom}_{B}(N, J^{i})\in {\rm add}DN$$
we see that $$({\rm Hom}_{A}(M, I^{t})\otimes_{k}{\rm End}_{B}(N))\oplus ({\rm End}_{A}(M)\otimes_{k}{\rm Hom}_{B}(N,J^{t}))$$ $$\oplus(\oplus_{i=0}^{t-1}({\rm Hom}_{A}(M, I^{i})\otimes_{k}{\rm Hom}_{B}(N, J^{t-1-i})))\in {\rm add}D(M\otimes_{k}N).$$ It follows that
$$D(M\otimes_{k}N)\mbox{-}{\rm domdim}{\rm End}_{A\otimes_{k}B}(M\otimes_{k}N)= {\rm min}\{DM\mbox{-}{\rm domdim}{\rm End}_{A}M,$$
$$DN\mbox{-}{\rm domdim}{\rm End}_{B}N\}.$$
Therefore, we get that
$${\rm Grigdim}(A\otimes_{k}B)\geq {\rm min}\{{\rm Grigdim}A,\ {\rm Grigdim}B\}.$$

\vskip 5pt

(3) \ Denote by $F: A\mbox{-}{\rm mod}\lxr B\mbox{-}{\rm mod}$ the Morita equivalence. We first claim that $F(M)$ is the Gorenstein
projective generator in $B\mbox{-}{\rm mod}$. Indeed, let
$$(P^{\bullet}, d)= \cdots \lxr P^{-1}\stackrel{d^{-1}}{\lxr} P^{0}\stackrel{d^{0}}{\lxr} P^{1}\stackrel{d^{1}}{\lxr} P^{2}\lxr \cdots$$ be a complete projective resolution of projective $A\mbox{-}$modules such that $M\cong {\rm Im}d^{0}$. Applying $F$, we get the exact sequence $F(P^{\bullet})$ of $B$-modules
$$\cdots \lxr F(P^{-1})\xrightarrow{F(d^{-1})} F(P^{0})\xrightarrow{F(d^{0})} F(P^{1})\xrightarrow{F(d^{1})} F(P^{2})\lxr \cdots$$
such that $F(M)\cong {\rm Im}F(d^{0})$, and moreover, there is the isomorphism
$${\rm Hom}_{B}(F(P^{\bullet}),\ Q)\cong {\rm Hom}_{A}(P^{\bullet},\ F^{-1}(Q))$$
for any projective $B\mbox{-}$module $Q$, which means that $F(P^{\bullet})$ is a complete projective resolution of projective $B\mbox{-}$modules. Hence, $F(M)$ is a Gorenstein projective $B$-module. Moreover, since $A\in {\rm add}M$, we get that $B\in {\rm add}F(M)$. Therefore, $F(M)$ is the Gorenstein projective generator in $B\mbox{-}{\rm mod}$.

\vskip 5pt

Next since $F$ is a Morita equivalence, there is the isomorphsim
$${\rm Ext}_{B}^{i}(F(M), F(M))\cong {\rm Ext}_{A}^{i}(M, M).$$
Thus we get that ${\rm Grigdim}A={\rm Grigdim}B$ from Lemma~\ref{domdim}.
\end{proof}
\end{prop}

\vskip 10pt

\begin{thm}\
\label{equivalences}
If $A\mbox{-}\underline{{\rm Gproj}}$ and $B\mbox{-}\underline{{\rm Gproj}}$ are stably equivalent as additive categories, then
$${\rm Grigdim}A= {\rm Grigdim}B.$$
Consequently,  if $A$ and $B$ are derived equivalent, then ${\rm Grigdim}A={\rm Grigdim}B$.
\begin{proof}\  Let $G: A\mbox{-}\underline{{\rm Gproj}}\lxr B\mbox{-}\underline{{\rm Gproj}}$ be a stable equivalence. Then we have that $G(M)\oplus B$ is the Gorenstein-projective generator in $B\mbox{-}{\rm mod}$. Since
$${\rm Ext}_{B}^{i}(G(M)\oplus B,\ G(M)\oplus B)\cong {\rm Ext}_{B}^{i}(G(M), G(M))\cong {\rm Ext}_{A}^{i}(M, M),$$ we get that ${\rm rd}_{A}M= {\rm rd}_{B}(G(M)\oplus B)$. Then  we have from Lemma~\ref{domdim} that $${\rm Grigdim}A= {\rm Grigdim}B.$$

\vskip 5pt

From \cite[Corollary 5.5]{HP}, we have that $A\mbox{-}\underline{{\rm Gproj}}$ and $B\mbox{-}\underline{{\rm Gproj}}$ are stably equivalent. Then ${\rm Grigdim}A={\rm Grigdim}B$ by the above arguments.
\end{proof}
\end{thm}

\vskip 10pt

\begin{exam}
\label{grid}
Let $A$ be a CM-finite non-CM-free Gorenstein algebra with the Gorenstein projective generator $G$. Then
$${\rm Grigdim}A\leq {\rm id}G+1.$$
\begin{proof}\ Let $m:={\rm id}G$. If $m=0$, then $G$ is injective, since $G$ is Gorenstein projective, there exists an exact sequence $0 \lxr G \lxr P \lxr C \lxr 0$ with $P$ projective. This means that $G$ is projective, and so $A$ is CM-free. Thus $m \geq 1$.

\vskip 5pt

When $m$ is infinite, there is nothing to show. Suppose
$m$ is finite. Then ${\rm Ext}^{m}_{A}(G, G)\neq 0$ and so ${\rm rd}_{A}G\leq m-1$.
Thus we have from above arguments that
$${\rm Grigdim}A\leq {\rm id}G-1+2={\rm id}G+1.$$
This completes the proof.
\end{proof}
\end{exam}

\vskip 10pt

In the following, we characterize the Gorenstein rigidity dimension  of an
algebra $A$ in terms of recollements of $A\mbox{-}\underline{{\rm Gproj}}$.

\vskip 10pt

\begin{thm}
\label{recollementgridity}
Assume that \[
\xymatrix@C=0.5cm{
B\mbox{-}{\rm \underline{Gproj}} \ar[rrr]^{\mathsf{i}} &&& A\mbox{-}{\rm \underline{Gproj}} \ar[rrr]^{\mathsf{e}}  \ar @/_1.5pc/[lll]_{\mathsf{q}}  \ar
 @/^1.5pc/[lll]^{\mathsf{p}} &&& C\mbox{-}{\rm \underline{Gproj}}
\ar @/_1.5pc/[lll]_{\mathsf{l}} \ar
 @/^1.5pc/[lll]^{\mathsf{r}}
 }\ \ \ \ \ \ \ (\ast)
\]
is a recollement of triangulated categories over the three algebras $A, B$ and $C$. Then the following hold.

\vskip 5pt

${\rm (1)}$ \ If $(\ast)$ admits the height 2 going upwards, then
$${\rm Grigdim}A\geq {\rm max} \{{\rm Grigdim}C,\  {\rm Grigdim}B\}.$$

\vskip 5pt

${\rm (2)}$ \ If $(\ast)$ admits the height 2 going downwards, then
$${\rm Grigdim}A\geq {\rm max} \{{\rm Grigdim}C, \ {\rm Grigdim}B\}.$$
\begin{proof}
We prove (1). Let $G$ be a Gorenstein projective generator in $C$-mod such that ${\rm rd}_{C}G=n-2$. Since $\mathsf{l}$ is fully faithful and exact, there are the isomorphisms
$${\rm Ext}^{i}_{A}(\mathsf{l}(G), \mathsf{l}(G))\cong {\rm Ext}_{C}^{i}(G, \mathsf{el}(G)) \cong {\rm Ext}_{C}^{i}(G, G)=0$$
for all $1\leq i\leq n-2$. This implies that ${\rm l}(G)\oplus A$ is a Gorenstein projective generator in A-mod such that $\rm{rd}_{A}({\rm l}(G)\oplus A)\geq n-2$, and therefore $\rm{Grigdim}A\geq n$.

\vskip 5pt

Let $E$ be a Gorenstein projective generator in $B\mbox{-}{\rm mod}$ such that ${\rm rd}_{B}E \geq n-2$. Since $\mathsf{i}$ is fully faithful and exact, there are isomorphisms $${\rm Ext}_{A}^{i}(\mathsf{i}(E), \mathsf{i}(E)) \cong {\rm Ext}^{i}_{B}(\mathsf{q}\mathsf{i}(E), E) \cong {\rm Ext}^{i}_{B}(E, E)=0$$ for all $1\leq i\leq n-2$. This implies that $\mathsf{i}(E)\oplus A$ is a Gorenstein projective generator in $A\mbox{-}{\rm mod}$ such that ${\rm rd}_{A}(\mathsf{i}(E)\oplus A) \geq n-2$ and therefore ${\rm Grigdim}A \geq n-2$.

\vskip 10pt

(2) is similarly proved.
\end{proof}
\end{thm}

\vskip 10pt

\begin{exam}
Let $A$ be a Gorenstein algebra and let $e$ and $f$ be two idempotent elements of $A$ such that $fAe=0$. Let $N:= Ae\otimes_{A}fA$ and
$\Lambda_{(0, 0)}=
\begin{pmatrix}
A & N\\
N & A
\end{pmatrix}
$ be the Morita context algebra.
Then $\rm{Grigdim} \Lambda_{(0, 0)} \geq \rm{Grigdim}A$.
\begin{proof}
By \cite[Example 3.1]{GX}, there exists an unbounded ladder
\[
\xymatrix@C=0.5cm{
A\mbox{-}{\rm \underline{Gproj}} \ar[rrr]^{\mathsf{i}} &&& \Lambda_{(0, 0)}\mbox{-}{\rm \underline{Gproj} } \ar[rrr]^{\mathsf{e}}  \ar @/_1.5pc/[lll]_{\mathsf{q}}  \ar
 @/^1.5pc/[lll]^{\mathsf{p}} &&&  A\mbox{-}{\rm \underline{Gproj}}.
\ar @/_1.5pc/[lll]_{\mathsf{l}} \ar
 @/^1.5pc/[lll]^{\mathsf{r}}
 }
 \]
Then the result immediately follows from Theorem~\ref{recollementgridity}.
\end{proof}
\end{exam}

\vskip 10pt

\begin{exam}
Let $A$ be a Gorenstein algebra. Let $\lambda: A \rightarrow B$ be a homological ring epimorphism which induces a recollement:
\[
\xymatrix@C=0.5cm{
D(B\mbox{-}{\rm Mod}) \ar[rrr]^{i_{\ast}} &&& D(A\mbox{-}{\rm Mod}) \ar[rrr]^{j^{\ast}}  \ar @/_1.5pc/[lll]_{i^{\ast}}  \ar
 @/^1.5pc/[lll]^{i^{!}} &&& D(C\mbox{-}{\rm Mod}),
\ar @/_1.5pc/[lll]_{j_{!}} \ar
 @/^1.5pc/[lll]^{j_{\ast}}
 }
\] such that $j_{!}$ restricts to $D^{b}(C\mbox{-}{\rm Mod})$. If ${\rm pd}_{A}B <\infty$, then
$${\rm Grigdim}A\geq {\rm max} \{{\rm Grigdim}C, \  {\rm Grigdim}B\}.$$
\end{exam}
\begin{proof}
By \cite[Theorem 3.1]{GX}, there is an unbounded ladder
\[
\xymatrix@C=0.5cm{
B\mbox{-}{\rm \underline{Gproj}} \ar[rrr]^{\mathsf{i}} &&& A\mbox{-}{\rm \underline{Gproj} } \ar[rrr]^{\mathsf{e}}  \ar @/_1.5pc/[lll]_{\mathsf{q}}  \ar
 @/^1.5pc/[lll]^{\mathsf{p}} &&&  C\mbox{-}{\rm \underline{Gproj}}.
\ar @/_1.5pc/[lll]_{\mathsf{l}} \ar
 @/^1.5pc/[lll]^{\mathsf{r}}
 }
\]
Then the result follows from  Theorem~\ref{recollementgridity}.
\end{proof}

\end{document}

%% file: rgb.tex
\definecolor{AliceBlue}{rgb}{0.94,0.97,1.00}
\definecolor{AntiqueWhite1}{rgb}{1.00,0.94,0.86}
\definecolor{AntiqueWhite2}{rgb}{0.93,0.87,0.80}
\definecolor{AntiqueWhite3}{rgb}{0.80,0.75,0.69}
\definecolor{AntiqueWhite4}{rgb}{0.55,0.51,0.47}
\definecolor{AntiqueWhite}{rgb}{0.98,0.92,0.84}
\definecolor{BlanchedAlmond}{rgb}{1.00,0.92,0.80}
\definecolor{BlueViolet}{rgb}{0.54,0.17,0.89}
\definecolor{CadetBlue1}{rgb}{0.60,0.96,1.00}
\definecolor{CadetBlue2}{rgb}{0.56,0.90,0.93}
\definecolor{CadetBlue3}{rgb}{0.48,0.77,0.80}
\definecolor{CadetBlue4}{rgb}{0.33,0.53,0.55}
\definecolor{CadetBlue}{rgb}{0.37,0.62,0.63}
\definecolor{CornflowerBlue}{rgb}{0.39,0.58,0.93}
\definecolor{DarkBlue}{rgb}{0.00,0.00,0.55}
\definecolor{DarkCyan}{rgb}{0.00,0.55,0.55}
\definecolor{DarkGoldenrod1}{rgb}{1.00,0.73,0.06}
\definecolor{DarkGoldenrod2}{rgb}{0.93,0.68,0.05}
\definecolor{DarkGoldenrod3}{rgb}{0.80,0.58,0.05}
\definecolor{DarkGoldenrod4}{rgb}{0.55,0.40,0.03}
\definecolor{DarkGoldenrod}{rgb}{0.72,0.53,0.04}
\definecolor{DarkGray}{rgb}{0.66,0.66,0.66}
\definecolor{DarkGreen}{rgb}{0.00,0.39,0.00}
\definecolor{DarkGrey}{rgb}{0.66,0.66,0.66}
\definecolor{DarkKhaki}{rgb}{0.74,0.72,0.42}
\definecolor{DarkMagenta}{rgb}{0.55,0.00,0.55}
\definecolor{DarkOliveGreen1}{rgb}{0.79,1.00,0.44}
\definecolor{DarkOliveGreen2}{rgb}{0.74,0.93,0.41}
\definecolor{DarkOliveGreen3}{rgb}{0.64,0.80,0.35}
\definecolor{DarkOliveGreen4}{rgb}{0.43,0.55,0.24}
\definecolor{DarkOliveGreen}{rgb}{0.33,0.42,0.18}
\definecolor{DarkOrange1}{rgb}{1.00,0.50,0.00}
\definecolor{DarkOrange2}{rgb}{0.93,0.46,0.00}
\definecolor{DarkOrange3}{rgb}{0.80,0.40,0.00}
\definecolor{DarkOrange4}{rgb}{0.55,0.27,0.00}
\definecolor{DarkOrange}{rgb}{1.00,0.55,0.00}
\definecolor{DarkOrchid1}{rgb}{0.75,0.24,1.00}
\definecolor{DarkOrchid2}{rgb}{0.70,0.23,0.93}
\definecolor{DarkOrchid3}{rgb}{0.60,0.20,0.80}
\definecolor{DarkOrchid4}{rgb}{0.41,0.13,0.55}
\definecolor{DarkOrchid}{rgb}{0.60,0.20,0.80}
\definecolor{DarkRed}{rgb}{0.55,0.00,0.00}
\definecolor{DarkSalmon}{rgb}{0.91,0.59,0.48}
\definecolor{DarkSeaGreen1}{rgb}{0.76,1.00,0.76}
\definecolor{DarkSeaGreen2}{rgb}{0.71,0.93,0.71}
\definecolor{DarkSeaGreen3}{rgb}{0.61,0.80,0.61}
\definecolor{DarkSeaGreen4}{rgb}{0.41,0.55,0.41}
\definecolor{DarkSeaGreen}{rgb}{0.56,0.74,0.56}
\definecolor{DarkSlateBlue}{rgb}{0.28,0.24,0.55}
\definecolor{DarkSlateGray1}{rgb}{0.59,1.00,1.00}
\definecolor{DarkSlateGray2}{rgb}{0.55,0.93,0.93}
\definecolor{DarkSlateGray3}{rgb}{0.47,0.80,0.80}
\definecolor{DarkSlateGray4}{rgb}{0.32,0.55,0.55}
\definecolor{DarkSlateGray}{rgb}{0.18,0.31,0.31}
\definecolor{DarkSlateGrey}{rgb}{0.18,0.31,0.31}
\definecolor{DarkTurquoise}{rgb}{0.00,0.81,0.82}
\definecolor{DarkViolet}{rgb}{0.58,0.00,0.83}
\definecolor{DeepPink1}{rgb}{1.00,0.08,0.58}
\definecolor{DeepPink2}{rgb}{0.93,0.07,0.54}
\definecolor{DeepPink3}{rgb}{0.80,0.06,0.46}
\definecolor{DeepPink4}{rgb}{0.55,0.04,0.31}
\definecolor{DeepPink}{rgb}{1.00,0.08,0.58}
\definecolor{DeepSkyBlue1}{rgb}{0.00,0.75,1.00}
\definecolor{DeepSkyBlue2}{rgb}{0.00,0.70,0.93}
\definecolor{DeepSkyBlue3}{rgb}{0.00,0.60,0.80}
\definecolor{DeepSkyBlue4}{rgb}{0.00,0.41,0.55}
\definecolor{DeepSkyBlue}{rgb}{0.00,0.75,1.00}
\definecolor{DimGray}{rgb}{0.41,0.41,0.41}
\definecolor{DimGrey}{rgb}{0.41,0.41,0.41}
\definecolor{DodgerBlue1}{rgb}{0.12,0.56,1.00}
\definecolor{DodgerBlue2}{rgb}{0.11,0.53,0.93}
\definecolor{DodgerBlue3}{rgb}{0.09,0.45,0.80}
\definecolor{DodgerBlue4}{rgb}{0.06,0.31,0.55}
\definecolor{DodgerBlue}{rgb}{0.12,0.56,1.00}
\definecolor{FloralWhite}{rgb}{1.00,0.98,0.94}
\definecolor{ForestGreen}{rgb}{0.13,0.55,0.13}
\definecolor{GhostWhite}{rgb}{0.97,0.97,1.00}
\definecolor{GreenYellow}{rgb}{0.68,1.00,0.18}
\definecolor{HotPink1}{rgb}{1.00,0.43,0.71}
\definecolor{HotPink2}{rgb}{0.93,0.42,0.65}
\definecolor{HotPink3}{rgb}{0.80,0.38,0.56}
\definecolor{HotPink4}{rgb}{0.55,0.23,0.38}
\definecolor{HotPink}{rgb}{1.00,0.41,0.71}
\definecolor{IndianRed1}{rgb}{1.00,0.42,0.42}
\definecolor{IndianRed2}{rgb}{0.93,0.39,0.39}
\definecolor{IndianRed3}{rgb}{0.80,0.33,0.33}
\definecolor{IndianRed4}{rgb}{0.55,0.23,0.23}
\definecolor{IndianRed}{rgb}{0.80,0.36,0.36}
\definecolor{LavenderBlush1}{rgb}{1.00,0.94,0.96}
\definecolor{LavenderBlush2}{rgb}{0.93,0.88,0.90}
\definecolor{LavenderBlush3}{rgb}{0.80,0.76,0.77}
\definecolor{LavenderBlush4}{rgb}{0.55,0.51,0.53}
\definecolor{LavenderBlush}{rgb}{1.00,0.94,0.96}
\definecolor{LawnGreen}{rgb}{0.49,0.99,0.00}
\definecolor{LemonChiffon1}{rgb}{1.00,0.98,0.80}
\definecolor{LemonChiffon2}{rgb}{0.93,0.91,0.75}
\definecolor{LemonChiffon3}{rgb}{0.80,0.79,0.65}
\definecolor{LemonChiffon4}{rgb}{0.55,0.54,0.44}
\definecolor{LemonChiffon}{rgb}{1.00,0.98,0.80}
\definecolor{LightBlue1}{rgb}{0.75,0.94,1.00}
\definecolor{LightBlue2}{rgb}{0.70,0.87,0.93}
\definecolor{LightBlue3}{rgb}{0.60,0.75,0.80}
\definecolor{LightBlue4}{rgb}{0.41,0.51,0.55}
\definecolor{LightBlue}{rgb}{0.68,0.85,0.90}
\definecolor{LightCoral}{rgb}{0.94,0.50,0.50}
\definecolor{LightCyan1}{rgb}{0.88,1.00,1.00}
\definecolor{LightCyan2}{rgb}{0.82,0.93,0.93}
\definecolor{LightCyan3}{rgb}{0.71,0.80,0.80}
\definecolor{LightCyan4}{rgb}{0.48,0.55,0.55}
\definecolor{LightCyan}{rgb}{0.88,1.00,1.00}
\definecolor{LightGoldenrod1}{rgb}{1.00,0.93,0.55}
\definecolor{LightGoldenrod2}{rgb}{0.93,0.86,0.51}
\definecolor{LightGoldenrod3}{rgb}{0.80,0.75,0.44}
\definecolor{LightGoldenrod4}{rgb}{0.55,0.51,0.30}
\definecolor{LightGoldenrodYellow}{rgb}{0.98,0.98,0.82}
\definecolor{LightGoldenrod}{rgb}{0.93,0.87,0.51}
\definecolor{LightGray}{rgb}{0.83,0.83,0.83}
\definecolor{LightGreen}{rgb}{0.56,0.93,0.56}
\definecolor{LightGrey}{rgb}{0.83,0.83,0.83}
\definecolor{LightPink1}{rgb}{1.00,0.68,0.73}
\definecolor{LightPink2}{rgb}{0.93,0.64,0.68}
\definecolor{LightPink3}{rgb}{0.80,0.55,0.58}
\definecolor{LightPink4}{rgb}{0.55,0.37,0.40}
\definecolor{LightPink}{rgb}{1.00,0.71,0.76}
\definecolor{LightSalmon1}{rgb}{1.00,0.63,0.48}
\definecolor{LightSalmon2}{rgb}{0.93,0.58,0.45}
\definecolor{LightSalmon3}{rgb}{0.80,0.51,0.38}
\definecolor{LightSalmon4}{rgb}{0.55,0.34,0.26}
\definecolor{LightSalmon}{rgb}{1.00,0.63,0.48}
\definecolor{LightSeaGreen}{rgb}{0.13,0.70,0.67}
\definecolor{LightSkyBlue1}{rgb}{0.69,0.89,1.00}
\definecolor{LightSkyBlue2}{rgb}{0.64,0.83,0.93}
\definecolor{LightSkyBlue3}{rgb}{0.55,0.71,0.80}
\definecolor{LightSkyBlue4}{rgb}{0.38,0.48,0.55}
\definecolor{LightSkyBlue}{rgb}{0.53,0.81,0.98}
\definecolor{LightSlateBlue}{rgb}{0.52,0.44,1.00}
\definecolor{LightSlateGray}{rgb}{0.47,0.53,0.60}
\definecolor{LightSlateGrey}{rgb}{0.47,0.53,0.60}
\definecolor{LightSteelBlue1}{rgb}{0.79,0.88,1.00}
\definecolor{LightSteelBlue2}{rgb}{0.74,0.82,0.93}
\definecolor{LightSteelBlue3}{rgb}{0.64,0.71,0.80}
\definecolor{LightSteelBlue4}{rgb}{0.43,0.48,0.55}
\definecolor{LightSteelBlue}{rgb}{0.69,0.77,0.87}
\definecolor{LightYellow1}{rgb}{1.00,1.00,0.88}
\definecolor{LightYellow2}{rgb}{0.93,0.93,0.82}
\definecolor{LightYellow3}{rgb}{0.80,0.80,0.71}
\definecolor{LightYellow4}{rgb}{0.55,0.55,0.48}
\definecolor{LightYellow}{rgb}{1.00,1.00,0.88}
\definecolor{LimeGreen}{rgb}{0.20,0.80,0.20}
\definecolor{MediumAquamarine}{rgb}{0.40,0.80,0.67}
\definecolor{MediumBlue}{rgb}{0.00,0.00,0.80}
\definecolor{MediumOrchid1}{rgb}{0.88,0.40,1.00}
\definecolor{MediumOrchid2}{rgb}{0.82,0.37,0.93}
\definecolor{MediumOrchid3}{rgb}{0.71,0.32,0.80}
\definecolor{MediumOrchid4}{rgb}{0.48,0.22,0.55}
\definecolor{MediumOrchid}{rgb}{0.73,0.33,0.83}
\definecolor{MediumPurple1}{rgb}{0.67,0.51,1.00}
\definecolor{MediumPurple2}{rgb}{0.62,0.47,0.93}
\definecolor{MediumPurple3}{rgb}{0.54,0.41,0.80}
\definecolor{MediumPurple4}{rgb}{0.36,0.28,0.55}
\definecolor{MediumPurple}{rgb}{0.58,0.44,0.86}
\definecolor{MediumSeaGreen}{rgb}{0.24,0.70,0.44}
\definecolor{MediumSlateBlue}{rgb}{0.48,0.41,0.93}
\definecolor{MediumSpringGreen}{rgb}{0.00,0.98,0.60}
\definecolor{MediumTurquoise}{rgb}{0.28,0.82,0.80}
\definecolor{MediumVioletRed}{rgb}{0.78,0.08,0.52}
\definecolor{MidnightBlue}{rgb}{0.10,0.10,0.44}
\definecolor{MintCream}{rgb}{0.96,1.00,0.98}
\definecolor{MistyRose1}{rgb}{1.00,0.89,0.88}
\definecolor{MistyRose2}{rgb}{0.93,0.84,0.82}
\definecolor{MistyRose3}{rgb}{0.80,0.72,0.71}
\definecolor{MistyRose4}{rgb}{0.55,0.49,0.48}
\definecolor{MistyRose}{rgb}{1.00,0.89,0.88}
\definecolor{NavajoWhite1}{rgb}{1.00,0.87,0.68}
\definecolor{NavajoWhite2}{rgb}{0.93,0.81,0.63}
\definecolor{NavajoWhite3}{rgb}{0.80,0.70,0.55}
\definecolor{NavajoWhite4}{rgb}{0.55,0.47,0.37}
\definecolor{NavajoWhite}{rgb}{1.00,0.87,0.68}
\definecolor{NavyBlue}{rgb}{0.00,0.00,0.50}
\definecolor{OldLace}{rgb}{0.99,0.96,0.90}
\definecolor{OliveDrab1}{rgb}{0.75,1.00,0.24}
\definecolor{OliveDrab2}{rgb}{0.70,0.93,0.23}
\definecolor{OliveDrab3}{rgb}{0.60,0.80,0.20}
\definecolor{OliveDrab4}{rgb}{0.41,0.55,0.13}
\definecolor{OliveDrab}{rgb}{0.42,0.56,0.14}
\definecolor{OrangeRed1}{rgb}{1.00,0.27,0.00}
\definecolor{OrangeRed2}{rgb}{0.93,0.25,0.00}
\definecolor{OrangeRed3}{rgb}{0.80,0.22,0.00}
\definecolor{OrangeRed4}{rgb}{0.55,0.15,0.00}
\definecolor{OrangeRed}{rgb}{1.00,0.27,0.00}
\definecolor{PaleGoldenrod}{rgb}{0.93,0.91,0.67}
\definecolor{PaleGreen1}{rgb}{0.60,1.00,0.60}
\definecolor{PaleGreen2}{rgb}{0.56,0.93,0.56}
\definecolor{PaleGreen3}{rgb}{0.49,0.80,0.49}
\definecolor{PaleGreen4}{rgb}{0.33,0.55,0.33}
\definecolor{PaleGreen}{rgb}{0.60,0.98,0.60}
\definecolor{PaleTurquoise1}{rgb}{0.73,1.00,1.00}
\definecolor{PaleTurquoise2}{rgb}{0.68,0.93,0.93}
\definecolor{PaleTurquoise3}{rgb}{0.59,0.80,0.80}
\definecolor{PaleTurquoise4}{rgb}{0.40,0.55,0.55}
\definecolor{PaleTurquoise}{rgb}{0.69,0.93,0.93}
\definecolor{PaleVioletRed1}{rgb}{1.00,0.51,0.67}
\definecolor{PaleVioletRed2}{rgb}{0.93,0.47,0.62}
\definecolor{PaleVioletRed3}{rgb}{0.80,0.41,0.54}
\definecolor{PaleVioletRed4}{rgb}{0.55,0.28,0.36}
\definecolor{PaleVioletRed}{rgb}{0.86,0.44,0.58}
\definecolor{PapayaWhip}{rgb}{1.00,0.94,0.84}
\definecolor{PeachPuff1}{rgb}{1.00,0.85,0.73}
\definecolor{PeachPuff2}{rgb}{0.93,0.80,0.68}
\definecolor{PeachPuff3}{rgb}{0.80,0.69,0.58}
\definecolor{PeachPuff4}{rgb}{0.55,0.47,0.40}
\definecolor{PeachPuff}{rgb}{1.00,0.85,0.73}
\definecolor{PowderBlue}{rgb}{0.69,0.88,0.90}
\definecolor{RosyBrown1}{rgb}{1.00,0.76,0.76}
\definecolor{RosyBrown2}{rgb}{0.93,0.71,0.71}
\definecolor{RosyBrown3}{rgb}{0.80,0.61,0.61}
\definecolor{RosyBrown4}{rgb}{0.55,0.41,0.41}
\definecolor{RosyBrown}{rgb}{0.74,0.56,0.56}
\definecolor{RoyalBlue1}{rgb}{0.28,0.46,1.00}
\definecolor{RoyalBlue2}{rgb}{0.26,0.43,0.93}
\definecolor{RoyalBlue3}{rgb}{0.23,0.37,0.80}
\definecolor{RoyalBlue4}{rgb}{0.15,0.25,0.55}
\definecolor{RoyalBlue}{rgb}{0.25,0.41,0.88}
\definecolor{SaddleBrown}{rgb}{0.55,0.27,0.07}
\definecolor{SandyBrown}{rgb}{0.96,0.64,0.38}
\definecolor{SeaGreen1}{rgb}{0.33,1.00,0.62}
\definecolor{SeaGreen2}{rgb}{0.31,0.93,0.58}
\definecolor{SeaGreen3}{rgb}{0.26,0.80,0.50}
\definecolor{SeaGreen4}{rgb}{0.18,0.55,0.34}
\definecolor{SeaGreen}{rgb}{0.18,0.55,0.34}
\definecolor{SkyBlue1}{rgb}{0.53,0.81,1.00}
\definecolor{SkyBlue2}{rgb}{0.49,0.75,0.93}
\definecolor{SkyBlue3}{rgb}{0.42,0.65,0.80}
\definecolor{SkyBlue4}{rgb}{0.29,0.44,0.55}
\definecolor{SkyBlue}{rgb}{0.53,0.81,0.92}
\definecolor{SlateBlue1}{rgb}{0.51,0.44,1.00}
\definecolor{SlateBlue2}{rgb}{0.48,0.40,0.93}
\definecolor{SlateBlue3}{rgb}{0.41,0.35,0.80}
\definecolor{SlateBlue4}{rgb}{0.28,0.24,0.55}
\definecolor{SlateBlue}{rgb}{0.42,0.35,0.80}
\definecolor{SlateGray1}{rgb}{0.78,0.89,1.00}
\definecolor{SlateGray2}{rgb}{0.73,0.83,0.93}
\definecolor{SlateGray3}{rgb}{0.62,0.71,0.80}
\definecolor{SlateGray4}{rgb}{0.42,0.48,0.55}
\definecolor{SlateGray}{rgb}{0.44,0.50,0.56}
\definecolor{SlateGrey}{rgb}{0.44,0.50,0.56}
\definecolor{SpringGreen1}{rgb}{0.00,1.00,0.50}
\definecolor{SpringGreen2}{rgb}{0.00,0.93,0.46}
\definecolor{SpringGreen3}{rgb}{0.00,0.80,0.40}
\definecolor{SpringGreen4}{rgb}{0.00,0.55,0.27}
\definecolor{SpringGreen}{rgb}{0.00,1.00,0.50}
\definecolor{SteelBlue1}{rgb}{0.39,0.72,1.00}
\definecolor{SteelBlue2}{rgb}{0.36,0.67,0.93}
\definecolor{SteelBlue3}{rgb}{0.31,0.58,0.80}
\definecolor{SteelBlue4}{rgb}{0.21,0.39,0.55}
\definecolor{SteelBlue}{rgb}{0.27,0.51,0.71}
\definecolor{VioletRed1}{rgb}{1.00,0.24,0.59}
\definecolor{VioletRed2}{rgb}{0.93,0.23,0.55}
\definecolor{VioletRed3}{rgb}{0.80,0.20,0.47}
\definecolor{VioletRed4}{rgb}{0.55,0.13,0.32}
\definecolor{VioletRed}{rgb}{0.82,0.13,0.56}
\definecolor{WhiteSmoke}{rgb}{0.96,0.96,0.96}
\definecolor{YellowGreen}{rgb}{0.60,0.80,0.20}
\definecolor{aliceblue}{rgb}{0.94,0.97,1.00}
\definecolor{antiquewhite}{rgb}{0.98,0.92,0.84}
\definecolor{aquamarine1}{rgb}{0.50,1.00,0.83}
\definecolor{aquamarine2}{rgb}{0.46,0.93,0.78}
\definecolor{aquamarine3}{rgb}{0.40,0.80,0.67}
\definecolor{aquamarine4}{rgb}{0.27,0.55,0.45}
\definecolor{aquamarine}{rgb}{0.50,1.00,0.83}
\definecolor{azure1}{rgb}{0.94,1.00,1.00}
\definecolor{azure2}{rgb}{0.88,0.93,0.93}
\definecolor{azure3}{rgb}{0.76,0.80,0.80}
\definecolor{azure4}{rgb}{0.51,0.55,0.55}
\definecolor{azure}{rgb}{0.94,1.00,1.00}
\definecolor{beige}{rgb}{0.96,0.96,0.86}
\definecolor{bisque1}{rgb}{1.00,0.89,0.77}
\definecolor{bisque2}{rgb}{0.93,0.84,0.72}
\definecolor{bisque3}{rgb}{0.80,0.72,0.62}
\definecolor{bisque4}{rgb}{0.55,0.49,0.42}
\definecolor{bisque}{rgb}{1.00,0.89,0.77}
\definecolor{black}{rgb}{0.00,0.00,0.00}
\definecolor{blanchedalmond}{rgb}{1.00,0.92,0.80}
\definecolor{blue1}{rgb}{0.00,0.00,1.00}
\definecolor{blue2}{rgb}{0.00,0.00,0.93}
\definecolor{blue3}{rgb}{0.00,0.00,0.80}
\definecolor{blue4}{rgb}{0.00,0.00,0.55}
\definecolor{blueviolet}{rgb}{0.54,0.17,0.89}
\definecolor{blue}{rgb}{0.00,0.00,1.00}
\definecolor{brown1}{rgb}{1.00,0.25,0.25}
\definecolor{brown2}{rgb}{0.93,0.23,0.23}
\definecolor{brown3}{rgb}{0.80,0.20,0.20}
\definecolor{brown4}{rgb}{0.55,0.14,0.14}
\definecolor{brown}{rgb}{0.65,0.16,0.16}
\definecolor{burlywood1}{rgb}{1.00,0.83,0.61}
\definecolor{burlywood2}{rgb}{0.93,0.77,0.57}
\definecolor{burlywood3}{rgb}{0.80,0.67,0.49}
\definecolor{burlywood4}{rgb}{0.55,0.45,0.33}
\definecolor{burlywood}{rgb}{0.87,0.72,0.53}
\definecolor{cadetblue}{rgb}{0.37,0.62,0.63}
\definecolor{chartreuse1}{rgb}{0.50,1.00,0.00}
\definecolor{chartreuse2}{rgb}{0.46,0.93,0.00}
\definecolor{chartreuse3}{rgb}{0.40,0.80,0.00}
\definecolor{chartreuse4}{rgb}{0.27,0.55,0.00}
\definecolor{chartreuse}{rgb}{0.50,1.00,0.00}
\definecolor{chocolate1}{rgb}{1.00,0.50,0.14}
\definecolor{chocolate2}{rgb}{0.93,0.46,0.13}
\definecolor{chocolate3}{rgb}{0.80,0.40,0.11}
\definecolor{chocolate4}{rgb}{0.55,0.27,0.07}
\definecolor{chocolate}{rgb}{0.82,0.41,0.12}
\definecolor{coral1}{rgb}{1.00,0.45,0.34}
\definecolor{coral2}{rgb}{0.93,0.42,0.31}
\definecolor{coral3}{rgb}{0.80,0.36,0.27}
\definecolor{coral4}{rgb}{0.55,0.24,0.18}
\definecolor{coral}{rgb}{1.00,0.50,0.31}
\definecolor{cornflowerblue}{rgb}{0.39,0.58,0.93}
\definecolor{cornsilk1}{rgb}{1.00,0.97,0.86}
\definecolor{cornsilk2}{rgb}{0.93,0.91,0.80}
\definecolor{cornsilk3}{rgb}{0.80,0.78,0.69}
\definecolor{cornsilk4}{rgb}{0.55,0.53,0.47}
\definecolor{cornsilk}{rgb}{1.00,0.97,0.86}
\definecolor{cyan1}{rgb}{0.00,1.00,1.00}
\definecolor{cyan2}{rgb}{0.00,0.93,0.93}
\definecolor{cyan3}{rgb}{0.00,0.80,0.80}
\definecolor{cyan4}{rgb}{0.00,0.55,0.55}
\definecolor{cyan}{rgb}{0.00,1.00,1.00}
\definecolor{darkblue}{rgb}{0.00,0.00,0.55}
\definecolor{darkcyan}{rgb}{0.00,0.55,0.55}
\definecolor{darkgoldenrod}{rgb}{0.72,0.53,0.04}
\definecolor{darkgray}{rgb}{0.66,0.66,0.66}
\definecolor{darkgreen}{rgb}{0.00,0.39,0.00}
\definecolor{darkgrey}{rgb}{0.66,0.66,0.66}
\definecolor{darkkhaki}{rgb}{0.74,0.72,0.42}
\definecolor{darkmagenta}{rgb}{0.55,0.00,0.55}
\definecolor{darkolive}{rgb}{0.33,0.42,0.18}
\definecolor{darkorange}{rgb}{1.00,0.55,0.00}
\definecolor{darkorchid}{rgb}{0.60,0.20,0.80}
\definecolor{darkred}{rgb}{0.55,0.00,0.00}
\definecolor{darksalmon}{rgb}{0.91,0.59,0.48}
\definecolor{darksea}{rgb}{0.56,0.74,0.56}
\definecolor{darkslate}{rgb}{0.18,0.31,0.31}
\definecolor{darkslate}{rgb}{0.18,0.31,0.31}
\definecolor{darkslate}{rgb}{0.28,0.24,0.55}
\definecolor{darkturquoise}{rgb}{0.00,0.81,0.82}
\definecolor{darkviolet}{rgb}{0.58,0.00,0.83}
\definecolor{deeppink}{rgb}{1.00,0.08,0.58}
\definecolor{deepsky}{rgb}{0.00,0.75,1.00}
\definecolor{dimgray}{rgb}{0.41,0.41,0.41}
\definecolor{dimgrey}{rgb}{0.41,0.41,0.41}
\definecolor{dodgerblue}{rgb}{0.12,0.56,1.00}
\definecolor{firebrick1}{rgb}{1.00,0.19,0.19}
\definecolor{firebrick2}{rgb}{0.93,0.17,0.17}
\definecolor{firebrick3}{rgb}{0.80,0.15,0.15}
\definecolor{firebrick4}{rgb}{0.55,0.10,0.10}
\definecolor{firebrick}{rgb}{0.70,0.13,0.13}
\definecolor{floralwhite}{rgb}{1.00,0.98,0.94}
\definecolor{forestgreen}{rgb}{0.13,0.55,0.13}
\definecolor{gainsboro}{rgb}{0.86,0.86,0.86}
\definecolor{ghostwhite}{rgb}{0.97,0.97,1.00}
\definecolor{gold1}{rgb}{1.00,0.84,0.00}
\definecolor{gold2}{rgb}{0.93,0.79,0.00}
\definecolor{gold3}{rgb}{0.80,0.68,0.00}
\definecolor{gold4}{rgb}{0.55,0.46,0.00}
\definecolor{goldenrod1}{rgb}{1.00,0.76,0.15}
\definecolor{goldenrod2}{rgb}{0.93,0.71,0.13}
\definecolor{goldenrod3}{rgb}{0.80,0.61,0.11}
\definecolor{goldenrod4}{rgb}{0.55,0.41,0.08}
\definecolor{goldenrod}{rgb}{0.85,0.65,0.13}
\definecolor{gold}{rgb}{1.00,0.84,0.00}
\definecolor{gray0}{rgb}{0.00,0.00,0.00}
\definecolor{gray100}{rgb}{1.00,1.00,1.00}
\definecolor{gray10}{rgb}{0.10,0.10,0.10}
\definecolor{gray11}{rgb}{0.11,0.11,0.11}
\definecolor{gray12}{rgb}{0.12,0.12,0.12}
\definecolor{gray13}{rgb}{0.13,0.13,0.13}
\definecolor{gray14}{rgb}{0.14,0.14,0.14}
\definecolor{gray15}{rgb}{0.15,0.15,0.15}
\definecolor{gray16}{rgb}{0.16,0.16,0.16}
\definecolor{gray17}{rgb}{0.17,0.17,0.17}
\definecolor{gray18}{rgb}{0.18,0.18,0.18}
\definecolor{gray19}{rgb}{0.19,0.19,0.19}
\definecolor{gray1}{rgb}{0.01,0.01,0.01}
\definecolor{gray20}{rgb}{0.20,0.20,0.20}
\definecolor{gray21}{rgb}{0.21,0.21,0.21}
\definecolor{gray22}{rgb}{0.22,0.22,0.22}
\definecolor{gray23}{rgb}{0.23,0.23,0.23}
\definecolor{gray24}{rgb}{0.24,0.24,0.24}
\definecolor{gray25}{rgb}{0.25,0.25,0.25}
\definecolor{gray26}{rgb}{0.26,0.26,0.26}
\definecolor{gray27}{rgb}{0.27,0.27,0.27}
\definecolor{gray28}{rgb}{0.28,0.28,0.28}
\definecolor{gray29}{rgb}{0.29,0.29,0.29}
\definecolor{gray2}{rgb}{0.02,0.02,0.02}
\definecolor{gray30}{rgb}{0.30,0.30,0.30}
\definecolor{gray31}{rgb}{0.31,0.31,0.31}
\definecolor{gray32}{rgb}{0.32,0.32,0.32}
\definecolor{gray33}{rgb}{0.33,0.33,0.33}
\definecolor{gray34}{rgb}{0.34,0.34,0.34}
\definecolor{gray35}{rgb}{0.35,0.35,0.35}
\definecolor{gray36}{rgb}{0.36,0.36,0.36}
\definecolor{gray37}{rgb}{0.37,0.37,0.37}
\definecolor{gray38}{rgb}{0.38,0.38,0.38}
\definecolor{gray39}{rgb}{0.39,0.39,0.39}
\definecolor{gray3}{rgb}{0.03,0.03,0.03}
\definecolor{gray40}{rgb}{0.40,0.40,0.40}
\definecolor{gray41}{rgb}{0.41,0.41,0.41}
\definecolor{gray42}{rgb}{0.42,0.42,0.42}
\definecolor{gray43}{rgb}{0.43,0.43,0.43}
\definecolor{gray44}{rgb}{0.44,0.44,0.44}
\definecolor{gray45}{rgb}{0.45,0.45,0.45}
\definecolor{gray46}{rgb}{0.46,0.46,0.46}
\definecolor{gray47}{rgb}{0.47,0.47,0.47}
\definecolor{gray48}{rgb}{0.48,0.48,0.48}
\definecolor{gray49}{rgb}{0.49,0.49,0.49}
\definecolor{gray4}{rgb}{0.04,0.04,0.04}
\definecolor{gray50}{rgb}{0.50,0.50,0.50}
\definecolor{gray51}{rgb}{0.51,0.51,0.51}
\definecolor{gray52}{rgb}{0.52,0.52,0.52}
\definecolor{gray53}{rgb}{0.53,0.53,0.53}
\definecolor{gray54}{rgb}{0.54,0.54,0.54}
\definecolor{gray55}{rgb}{0.55,0.55,0.55}
\definecolor{gray56}{rgb}{0.56,0.56,0.56}
\definecolor{gray57}{rgb}{0.57,0.57,0.57}
\definecolor{gray58}{rgb}{0.58,0.58,0.58}
\definecolor{gray59}{rgb}{0.59,0.59,0.59}
\definecolor{gray5}{rgb}{0.05,0.05,0.05}
\definecolor{gray60}{rgb}{0.60,0.60,0.60}
\definecolor{gray61}{rgb}{0.61,0.61,0.61}
\definecolor{gray62}{rgb}{0.62,0.62,0.62}
\definecolor{gray63}{rgb}{0.63,0.63,0.63}
\definecolor{gray64}{rgb}{0.64,0.64,0.64}
\definecolor{gray65}{rgb}{0.65,0.65,0.65}
\definecolor{gray66}{rgb}{0.66,0.66,0.66}
\definecolor{gray67}{rgb}{0.67,0.67,0.67}
\definecolor{gray68}{rgb}{0.68,0.68,0.68}
\definecolor{gray69}{rgb}{0.69,0.69,0.69}
\definecolor{gray6}{rgb}{0.06,0.06,0.06}
\definecolor{gray70}{rgb}{0.70,0.70,0.70}
\definecolor{gray71}{rgb}{0.71,0.71,0.71}
\definecolor{gray72}{rgb}{0.72,0.72,0.72}
\definecolor{gray73}{rgb}{0.73,0.73,0.73}
\definecolor{gray74}{rgb}{0.74,0.74,0.74}
\definecolor{gray75}{rgb}{0.75,0.75,0.75}
\definecolor{gray76}{rgb}{0.76,0.76,0.76}
\definecolor{gray77}{rgb}{0.77,0.77,0.77}
\definecolor{gray78}{rgb}{0.78,0.78,0.78}
\definecolor{gray79}{rgb}{0.79,0.79,0.79}
\definecolor{gray7}{rgb}{0.07,0.07,0.07}
\definecolor{gray80}{rgb}{0.80,0.80,0.80}
\definecolor{gray81}{rgb}{0.81,0.81,0.81}
\definecolor{gray82}{rgb}{0.82,0.82,0.82}
\definecolor{gray83}{rgb}{0.83,0.83,0.83}
\definecolor{gray84}{rgb}{0.84,0.84,0.84}
\definecolor{gray85}{rgb}{0.85,0.85,0.85}
\definecolor{gray86}{rgb}{0.86,0.86,0.86}
\definecolor{gray87}{rgb}{0.87,0.87,0.87}
\definecolor{gray88}{rgb}{0.88,0.88,0.88}
\definecolor{gray89}{rgb}{0.89,0.89,0.89}
\definecolor{gray8}{rgb}{0.08,0.08,0.08}
\definecolor{gray90}{rgb}{0.90,0.90,0.90}
\definecolor{gray91}{rgb}{0.91,0.91,0.91}
\definecolor{gray92}{rgb}{0.92,0.92,0.92}
\definecolor{gray93}{rgb}{0.93,0.93,0.93}
\definecolor{gray94}{rgb}{0.94,0.94,0.94}
\definecolor{gray95}{rgb}{0.95,0.95,0.95}
\definecolor{gray96}{rgb}{0.96,0.96,0.96}
\definecolor{gray97}{rgb}{0.97,0.97,0.97}
\definecolor{gray98}{rgb}{0.98,0.98,0.98}
\definecolor{gray99}{rgb}{0.99,0.99,0.99}
\definecolor{gray9}{rgb}{0.09,0.09,0.09}
\definecolor{gray}{rgb}{0.75,0.75,0.75}
\definecolor{green1}{rgb}{0.00,1.00,0.00}
\definecolor{green2}{rgb}{0.00,0.93,0.00}
\definecolor{green3}{rgb}{0.00,0.80,0.00}
\definecolor{green4}{rgb}{0.00,0.55,0.00}
\definecolor{greenyellow}{rgb}{0.68,1.00,0.18}
\definecolor{green}{rgb}{0.00,1.00,0.00}
\definecolor{grey0}{rgb}{0.00,0.00,0.00}
\definecolor{grey100}{rgb}{1.00,1.00,1.00}
\definecolor{grey10}{rgb}{0.10,0.10,0.10}
\definecolor{grey11}{rgb}{0.11,0.11,0.11}
\definecolor{grey12}{rgb}{0.12,0.12,0.12}
\definecolor{grey13}{rgb}{0.13,0.13,0.13}
\definecolor{grey14}{rgb}{0.14,0.14,0.14}
\definecolor{grey15}{rgb}{0.15,0.15,0.15}
\definecolor{grey16}{rgb}{0.16,0.16,0.16}
\definecolor{grey17}{rgb}{0.17,0.17,0.17}
\definecolor{grey18}{rgb}{0.18,0.18,0.18}
\definecolor{grey19}{rgb}{0.19,0.19,0.19}
\definecolor{grey1}{rgb}{0.01,0.01,0.01}
\definecolor{grey20}{rgb}{0.20,0.20,0.20}
\definecolor{grey21}{rgb}{0.21,0.21,0.21}
\definecolor{grey22}{rgb}{0.22,0.22,0.22}
\definecolor{grey23}{rgb}{0.23,0.23,0.23}
\definecolor{grey24}{rgb}{0.24,0.24,0.24}
\definecolor{grey25}{rgb}{0.25,0.25,0.25}
\definecolor{grey26}{rgb}{0.26,0.26,0.26}
\definecolor{grey27}{rgb}{0.27,0.27,0.27}
\definecolor{grey28}{rgb}{0.28,0.28,0.28}
\definecolor{grey29}{rgb}{0.29,0.29,0.29}
\definecolor{grey2}{rgb}{0.02,0.02,0.02}
\definecolor{grey30}{rgb}{0.30,0.30,0.30}
\definecolor{grey31}{rgb}{0.31,0.31,0.31}
\definecolor{grey32}{rgb}{0.32,0.32,0.32}
\definecolor{grey33}{rgb}{0.33,0.33,0.33}
\definecolor{grey34}{rgb}{0.34,0.34,0.34}
\definecolor{grey35}{rgb}{0.35,0.35,0.35}
\definecolor{grey36}{rgb}{0.36,0.36,0.36}
\definecolor{grey37}{rgb}{0.37,0.37,0.37}
\definecolor{grey38}{rgb}{0.38,0.38,0.38}
\definecolor{grey39}{rgb}{0.39,0.39,0.39}
\definecolor{grey3}{rgb}{0.03,0.03,0.03}
\definecolor{grey40}{rgb}{0.40,0.40,0.40}
\definecolor{grey41}{rgb}{0.41,0.41,0.41}
\definecolor{grey42}{rgb}{0.42,0.42,0.42}
\definecolor{grey43}{rgb}{0.43,0.43,0.43}
\definecolor{grey44}{rgb}{0.44,0.44,0.44}
\definecolor{grey45}{rgb}{0.45,0.45,0.45}
\definecolor{grey46}{rgb}{0.46,0.46,0.46}
\definecolor{grey47}{rgb}{0.47,0.47,0.47}
\definecolor{grey48}{rgb}{0.48,0.48,0.48}
\definecolor{grey49}{rgb}{0.49,0.49,0.49}
\definecolor{grey4}{rgb}{0.04,0.04,0.04}
\definecolor{grey50}{rgb}{0.50,0.50,0.50}
\definecolor{grey51}{rgb}{0.51,0.51,0.51}
\definecolor{grey52}{rgb}{0.52,0.52,0.52}
\definecolor{grey53}{rgb}{0.53,0.53,0.53}
\definecolor{grey54}{rgb}{0.54,0.54,0.54}
\definecolor{grey55}{rgb}{0.55,0.55,0.55}
\definecolor{grey56}{rgb}{0.56,0.56,0.56}
\definecolor{grey57}{rgb}{0.57,0.57,0.57}
\definecolor{grey58}{rgb}{0.58,0.58,0.58}
\definecolor{grey59}{rgb}{0.59,0.59,0.59}
\definecolor{grey5}{rgb}{0.05,0.05,0.05}
\definecolor{grey60}{rgb}{0.60,0.60,0.60}
\definecolor{grey61}{rgb}{0.61,0.61,0.61}
\definecolor{grey62}{rgb}{0.62,0.62,0.62}
\definecolor{grey63}{rgb}{0.63,0.63,0.63}
\definecolor{grey64}{rgb}{0.64,0.64,0.64}
\definecolor{grey65}{rgb}{0.65,0.65,0.65}
\definecolor{grey66}{rgb}{0.66,0.66,0.66}
\definecolor{grey67}{rgb}{0.67,0.67,0.67}
\definecolor{grey68}{rgb}{0.68,0.68,0.68}
\definecolor{grey69}{rgb}{0.69,0.69,0.69}
\definecolor{grey6}{rgb}{0.06,0.06,0.06}
\definecolor{grey70}{rgb}{0.70,0.70,0.70}
\definecolor{grey71}{rgb}{0.71,0.71,0.71}
\definecolor{grey72}{rgb}{0.72,0.72,0.72}
\definecolor{grey73}{rgb}{0.73,0.73,0.73}
\definecolor{grey74}{rgb}{0.74,0.74,0.74}
\definecolor{grey75}{rgb}{0.75,0.75,0.75}
\definecolor{grey76}{rgb}{0.76,0.76,0.76}
\definecolor{grey77}{rgb}{0.77,0.77,0.77}
\definecolor{grey78}{rgb}{0.78,0.78,0.78}
\definecolor{grey79}{rgb}{0.79,0.79,0.79}
\definecolor{grey7}{rgb}{0.07,0.07,0.07}
\definecolor{grey80}{rgb}{0.80,0.80,0.80}
\definecolor{grey81}{rgb}{0.81,0.81,0.81}
\definecolor{grey82}{rgb}{0.82,0.82,0.82}
\definecolor{grey83}{rgb}{0.83,0.83,0.83}
\definecolor{grey84}{rgb}{0.84,0.84,0.84}
\definecolor{grey85}{rgb}{0.85,0.85,0.85}
\definecolor{grey86}{rgb}{0.86,0.86,0.86}
\definecolor{grey87}{rgb}{0.87,0.87,0.87}
\definecolor{grey88}{rgb}{0.88,0.88,0.88}
\definecolor{grey89}{rgb}{0.89,0.89,0.89}
\definecolor{grey8}{rgb}{0.08,0.08,0.08}
\definecolor{grey90}{rgb}{0.90,0.90,0.90}
\definecolor{grey91}{rgb}{0.91,0.91,0.91}
\definecolor{grey92}{rgb}{0.92,0.92,0.92}
\definecolor{grey93}{rgb}{0.93,0.93,0.93}
\definecolor{grey94}{rgb}{0.94,0.94,0.94}
\definecolor{grey95}{rgb}{0.95,0.95,0.95}
\definecolor{grey96}{rgb}{0.96,0.96,0.96}
\definecolor{grey97}{rgb}{0.97,0.97,0.97}
\definecolor{grey98}{rgb}{0.98,0.98,0.98}
\definecolor{grey99}{rgb}{0.99,0.99,0.99}
\definecolor{grey9}{rgb}{0.09,0.09,0.09}
\definecolor{grey}{rgb}{0.75,0.75,0.75}
\definecolor{honeydew1}{rgb}{0.94,1.00,0.94}
\definecolor{honeydew2}{rgb}{0.88,0.93,0.88}
\definecolor{honeydew3}{rgb}{0.76,0.80,0.76}
\definecolor{honeydew4}{rgb}{0.51,0.55,0.51}
\definecolor{honeydew}{rgb}{0.94,1.00,0.94}
\definecolor{hotpink}{rgb}{1.00,0.41,0.71}
\definecolor{indianred}{rgb}{0.80,0.36,0.36}
\definecolor{ivory1}{rgb}{1.00,1.00,0.94}
\definecolor{ivory2}{rgb}{0.93,0.93,0.88}
\definecolor{ivory3}{rgb}{0.80,0.80,0.76}
\definecolor{ivory4}{rgb}{0.55,0.55,0.51}
\definecolor{ivory}{rgb}{1.00,1.00,0.94}
\definecolor{khaki1}{rgb}{1.00,0.96,0.56}
\definecolor{khaki2}{rgb}{0.93,0.90,0.52}
\definecolor{khaki3}{rgb}{0.80,0.78,0.45}
\definecolor{khaki4}{rgb}{0.55,0.53,0.31}
\definecolor{khaki}{rgb}{0.94,0.90,0.55}
\definecolor{lavenderblush}{rgb}{1.00,0.94,0.96}
\definecolor{lavender}{rgb}{0.90,0.90,0.98}
\definecolor{lawngreen}{rgb}{0.49,0.99,0.00}
\definecolor{lemonchiffon}{rgb}{1.00,0.98,0.80}
\definecolor{lightblue}{rgb}{0.68,0.85,0.90}
\definecolor{lightcoral}{rgb}{0.94,0.50,0.50}
\definecolor{lightcyan}{rgb}{0.88,1.00,1.00}
\definecolor{lightgoldenrod}{rgb}{0.93,0.87,0.51}
\definecolor{lightgoldenrod}{rgb}{0.98,0.98,0.82}
\definecolor{lightgray}{rgb}{0.83,0.83,0.83}
\definecolor{lightgreen}{rgb}{0.56,0.93,0.56}
\definecolor{lightgrey}{rgb}{0.83,0.83,0.83}
\definecolor{lightpink}{rgb}{1.00,0.71,0.76}
\definecolor{lightsalmon}{rgb}{1.00,0.63,0.48}
\definecolor{lightsea}{rgb}{0.13,0.70,0.67}
\definecolor{lightsky}{rgb}{0.53,0.81,0.98}
\definecolor{lightslate}{rgb}{0.47,0.53,0.60}
\definecolor{lightslate}{rgb}{0.47,0.53,0.60}
\definecolor{lightslate}{rgb}{0.52,0.44,1.00}
\definecolor{lightsteel}{rgb}{0.69,0.77,0.87}
\definecolor{lightyellow}{rgb}{1.00,1.00,0.88}
\definecolor{limegreen}{rgb}{0.20,0.80,0.20}
\definecolor{linen}{rgb}{0.98,0.94,0.90}
\definecolor{magenta1}{rgb}{1.00,0.00,1.00}
\definecolor{magenta2}{rgb}{0.93,0.00,0.93}
\definecolor{magenta3}{rgb}{0.80,0.00,0.80}
\definecolor{magenta4}{rgb}{0.55,0.00,0.55}
\definecolor{magenta}{rgb}{1.00,0.00,1.00}
\definecolor{maroon1}{rgb}{1.00,0.20,0.70}
\definecolor{maroon2}{rgb}{0.93,0.19,0.65}
\definecolor{maroon3}{rgb}{0.80,0.16,0.56}
\definecolor{maroon4}{rgb}{0.55,0.11,0.38}
\definecolor{maroon}{rgb}{0.69,0.19,0.38}
\definecolor{mediumaquamarine}{rgb}{0.40,0.80,0.67}
\definecolor{mediumblue}{rgb}{0.00,0.00,0.80}
\definecolor{mediumorchid}{rgb}{0.73,0.33,0.83}
\definecolor{mediumpurple}{rgb}{0.58,0.44,0.86}
\definecolor{mediumsea}{rgb}{0.24,0.70,0.44}
\definecolor{mediumslate}{rgb}{0.48,0.41,0.93}
\definecolor{mediumspring}{rgb}{0.00,0.98,0.60}
\definecolor{mediumturquoise}{rgb}{0.28,0.82,0.80}
\definecolor{mediumviolet}{rgb}{0.78,0.08,0.52}
\definecolor{midnightblue}{rgb}{0.10,0.10,0.44}
\definecolor{mintcream}{rgb}{0.96,1.00,0.98}
\definecolor{mistyrose}{rgb}{1.00,0.89,0.88}
\definecolor{moccasin}{rgb}{1.00,0.89,0.71}
\definecolor{navajowhite}{rgb}{1.00,0.87,0.68}
\definecolor{navyblue}{rgb}{0.00,0.00,0.50}
\definecolor{navy}{rgb}{0.00,0.00,0.50}
\definecolor{oldlace}{rgb}{0.99,0.96,0.90}
\definecolor{olivedrab}{rgb}{0.42,0.56,0.14}
\definecolor{orange1}{rgb}{1.00,0.65,0.00}
\definecolor{orange2}{rgb}{0.93,0.60,0.00}
\definecolor{orange3}{rgb}{0.80,0.52,0.00}
\definecolor{orange4}{rgb}{0.55,0.35,0.00}
\definecolor{orangered}{rgb}{1.00,0.27,0.00}
\definecolor{orange}{rgb}{1.00,0.65,0.00}
\definecolor{orchid1}{rgb}{1.00,0.51,0.98}
\definecolor{orchid2}{rgb}{0.93,0.48,0.91}
\definecolor{orchid3}{rgb}{0.80,0.41,0.79}
\definecolor{orchid4}{rgb}{0.55,0.28,0.54}
\definecolor{orchid}{rgb}{0.85,0.44,0.84}
\definecolor{palegoldenrod}{rgb}{0.93,0.91,0.67}
\definecolor{palegreen}{rgb}{0.60,0.98,0.60}
\definecolor{paleturquoise}{rgb}{0.69,0.93,0.93}
\definecolor{paleviolet}{rgb}{0.86,0.44,0.58}
\definecolor{papayawhip}{rgb}{1.00,0.94,0.84}
\definecolor{peachpuff}{rgb}{1.00,0.85,0.73}
\definecolor{peru}{rgb}{0.80,0.52,0.25}
\definecolor{pink1}{rgb}{1.00,0.71,0.77}
\definecolor{pink2}{rgb}{0.93,0.66,0.72}
\definecolor{pink3}{rgb}{0.80,0.57,0.62}
\definecolor{pink4}{rgb}{0.55,0.39,0.42}
\definecolor{pink}{rgb}{1.00,0.75,0.80}
\definecolor{plum1}{rgb}{1.00,0.73,1.00}
\definecolor{plum2}{rgb}{0.93,0.68,0.93}
\definecolor{plum3}{rgb}{0.80,0.59,0.80}
\definecolor{plum4}{rgb}{0.55,0.40,0.55}
\definecolor{plum}{rgb}{0.87,0.63,0.87}
\definecolor{powderblue}{rgb}{0.69,0.88,0.90}
\definecolor{purple1}{rgb}{0.61,0.19,1.00}
\definecolor{purple2}{rgb}{0.57,0.17,0.93}
\definecolor{purple3}{rgb}{0.49,0.15,0.80}
\definecolor{purple4}{rgb}{0.33,0.10,0.55}
\definecolor{purple}{rgb}{0.63,0.13,0.94}
\definecolor{red1}{rgb}{1.00,0.00,0.00}
\definecolor{red2}{rgb}{0.93,0.00,0.00}
\definecolor{red3}{rgb}{0.80,0.00,0.00}
\definecolor{red4}{rgb}{0.55,0.00,0.00}
\definecolor{red}{rgb}{1.00,0.00,0.00}
\definecolor{rosybrown}{rgb}{0.74,0.56,0.56}
\definecolor{royalblue}{rgb}{0.25,0.41,0.88}
\definecolor{saddlebrown}{rgb}{0.55,0.27,0.07}
\definecolor{salmon1}{rgb}{1.00,0.55,0.41}
\definecolor{salmon2}{rgb}{0.93,0.51,0.38}
\definecolor{salmon3}{rgb}{0.80,0.44,0.33}
\definecolor{salmon4}{rgb}{0.55,0.30,0.22}
\definecolor{salmon}{rgb}{0.98,0.50,0.45}
\definecolor{sandybrown}{rgb}{0.96,0.64,0.38}
\definecolor{seagreen}{rgb}{0.18,0.55,0.34}
\definecolor{seashell1}{rgb}{1.00,0.96,0.93}
\definecolor{seashell2}{rgb}{0.93,0.90,0.87}
\definecolor{seashell3}{rgb}{0.80,0.77,0.75}
\definecolor{seashell4}{rgb}{0.55,0.53,0.51}
\definecolor{seashell}{rgb}{1.00,0.96,0.93}
\definecolor{sienna1}{rgb}{1.00,0.51,0.28}
\definecolor{sienna2}{rgb}{0.93,0.47,0.26}
\definecolor{sienna3}{rgb}{0.80,0.41,0.22}
\definecolor{sienna4}{rgb}{0.55,0.28,0.15}
\definecolor{sienna}{rgb}{0.63,0.32,0.18}
\definecolor{skyblue}{rgb}{0.53,0.81,0.92}
\definecolor{slateblue}{rgb}{0.42,0.35,0.80}
\definecolor{slategray}{rgb}{0.44,0.50,0.56}
\definecolor{slategrey}{rgb}{0.44,0.50,0.56}
\definecolor{snow1}{rgb}{1.00,0.98,0.98}
\definecolor{snow2}{rgb}{0.93,0.91,0.91}
\definecolor{snow3}{rgb}{0.80,0.79,0.79}
\definecolor{snow4}{rgb}{0.55,0.54,0.54}
\definecolor{snow}{rgb}{1.00,0.98,0.98}
\definecolor{springgreen}{rgb}{0.00,1.00,0.50}
\definecolor{steelblue}{rgb}{0.27,0.51,0.71}
\definecolor{tan1}{rgb}{1.00,0.65,0.31}
\definecolor{tan2}{rgb}{0.93,0.60,0.29}
\definecolor{tan3}{rgb}{0.80,0.52,0.25}
\definecolor{tan4}{rgb}{0.55,0.35,0.17}
\definecolor{tan}{rgb}{0.82,0.71,0.55}
\definecolor{thistle1}{rgb}{1.00,0.88,1.00}
\definecolor{thistle2}{rgb}{0.93,0.82,0.93}
\definecolor{thistle3}{rgb}{0.80,0.71,0.80}
\definecolor{thistle4}{rgb}{0.55,0.48,0.55}
\definecolor{thistle}{rgb}{0.85,0.75,0.85}
\definecolor{tomato1}{rgb}{1.00,0.39,0.28}
\definecolor{tomato2}{rgb}{0.93,0.36,0.26}
\definecolor{tomato3}{rgb}{0.80,0.31,0.22}
\definecolor{tomato4}{rgb}{0.55,0.21,0.15}
\definecolor{tomato}{rgb}{1.00,0.39,0.28}
\definecolor{turquoise1}{rgb}{0.00,0.96,1.00}
\definecolor{turquoise2}{rgb}{0.00,0.90,0.93}
\definecolor{turquoise3}{rgb}{0.00,0.77,0.80}
\definecolor{turquoise4}{rgb}{0.00,0.53,0.55}
\definecolor{turquoise}{rgb}{0.25,0.88,0.82}
\definecolor{violetred}{rgb}{0.82,0.13,0.56}
\definecolor{violet}{rgb}{0.93,0.51,0.93}
\definecolor{wheat1}{rgb}{1.00,0.91,0.73}
\definecolor{wheat2}{rgb}{0.93,0.85,0.68}
\definecolor{wheat3}{rgb}{0.80,0.73,0.59}
\definecolor{wheat4}{rgb}{0.55,0.49,0.40}
\definecolor{wheat}{rgb}{0.96,0.87,0.70}
\definecolor{whitesmoke}{rgb}{0.96,0.96,0.96}
\definecolor{white}{rgb}{1.00,1.00,1.00}
\definecolor{yellow1}{rgb}{1.00,1.00,0.00}
\definecolor{yellow2}{rgb}{0.93,0.93,0.00}
\definecolor{yellow3}{rgb}{0.80,0.80,0.00}
\definecolor{yellow4}{rgb}{0.55,0.55,0.00}
\definecolor{yellowgreen}{rgb}{0.60,0.80,0.20}
\definecolor{yellow}{rgb}{1.00,1.00,0.00}